\documentclass{article}

\usepackage{amsmath,amssymb}
\usepackage[latin1]{inputenc}
\usepackage[dvips]{graphics}
\input{xy}
\xyoption{poly}
\xyoption{2cell}
\xyoption{all}

\newcommand{\ZZ}{\mathbb{Z}}

\newcommand{\qed}{$\ \hfill\Box $\bigskip}

\newcommand{\edges}{E'} 
\newcommand{\CD}{\mathcal{C}_{D_n}}
\newcommand{\CA}{\mathcal{C}_{A}}
\newcommand{\pos}{\textup{pos}\,}
\newcommand{\level}{\textup{level}\,}
\newcommand{\mesh}{\mathcal{M}}

\newcommand{\Hom}{\textup{Hom}\,}

\newcommand{\Ext}{\textup{Ext}}
\newcommand{\End}{\textup{End}\,}

\newcommand{\ind}{\textup{ind}\,}

\newcommand{\DA}{\mathcal{D}^b\! A}

\newcommand{\ssi}{\Leftrightarrow}

\newcommand{\za}{\alpha}
\newcommand{\zb}{\beta}
\newcommand{\zd}{\delta}
\newcommand{\ze}{\epsilon}
\newcommand{\zg}{\gamma}
\newcommand{\zG}{\Gamma}
\newcommand{\zl}{\lambda}
\newcommand{\zs}{\sigma}

\newtheorem{thm}{Theorem}[section]
\newtheorem{prop}[thm]{Proposition}

\newtheorem{cor}[thm]{Corollary}
\newtheorem{lem}[thm]{Lemma}
\newtheorem{rem}[thm]{Remark}

\newenvironment{pf}{{Proof}.}

\begin{document}
\title{A geometric model for cluster categories of type $D_n$}
\author{R. Schiffler}
\date{}
\maketitle


\begin{abstract}
We give a geometric realization of cluster categories of type $D_n$
using a polygon with $n$ vertices and one puncture in its center as a
model. In this realization, the indecomposable objects of the cluster
category correspond to certain homotopy classes of paths between two
vertices. 
 \end{abstract}
\maketitle

\setcounter{section}{0}
\addtocounter{section}{-1}


\begin{section}{Introduction}\label{section intro}
Cluster categories were introduced in \cite{BMRRT} and, independently,
in \cite{CCS} for type $A_n$,  as a means for better understanding of
the cluster algebras of Fomin and Zelevinsky \cite{FZ1,FZ2}. Since
then cluster categories have been the subject of many investigations,
see, for instance,
\cite{ABST1,ABST2,BMR1,BMR2,BMRT,CC,CCS2,CK1,CK2,K,KZ,Zhu1}. 

In the approach of \cite{BMRRT}, the cluster category $\CA$ is defined
as the quotient $\DA/F$ of the derived category $\DA$ of a hereditary
algebra $A$ by the endofunctor $F=\tau^{-1}_{\DA}\,[1]$, where
$\tau_{\DA}$ is the Auslander-Reiten translation and $[1]$ is the
shift. 
On the other hand, in the approach of \cite{CCS}, which is only valid
in type $A_n$, the cluster category is realized by an ad-hoc method 
as a category of
diagonals of a regular polygon with $n+3$ vertices. The morphisms
between diagonals are constructed geometrically using so called
elementary moves and mesh relations. 
In that realization, clusters are in one-to-one correspondence with
triangulations of the polygon and mutations are given by flips of 
diagonals in  the triangulation. 
Recently, Baur and Marsh \cite{Baur-Marsh} have generalized this model
to $m$-cluster categories of type $A_n$.

In this paper, we give a geometric realization of the cluster
categories of type $D_n$ in the spirit of \cite{CCS}. 
The polygon with $(n+3)$ vertices has to be replaced by a polygon with
$n$ vertices and one puncture in the center, and instead of looking at
diagonals, which are straight lines between two vertices, one has to
consider homotopy classes of paths between two
vertices, which we will call edges.
This  punctured polygon model has appeared
recently in the work of Fomin, Shapiro and Thurston \cite{FST} on the
relation between cluster algebras and triangulated surfaces. Let us
point out that they work in a vastly more general context and the
punctured polygon is only one example of their theory.
We define the cluster category by an ad-hoc method as the category of
(tagged) edges inside the punctured polygon. Morphisms are defined
using so-called elementary moves and mesh relations,  which are
generalizations of the 
elementary moves and mesh relations of \cite{CCS}. Our main results
are the equivalence of the category of tagged edges and the cluster
category of \cite{BMRRT}, see Theorem \ref{thm 1}, and the realization of the
dimension of $\Ext^1$ of tagged edges as the number of crossings
between  the same  tagged edges, see Theorem \ref{thm 2}.

The article is organized as follows. After a brief preliminary section,
in which we fix the notations and recall some concepts needed
later,
section \ref{sect 2} is devoted to the definition of the category
$\mathcal{C}$  of tagged
edges. In section \ref{sect 3}, we show the equivalence of this
category and the cluster category and in section \ref{sect 4}, we study
$\Ext^1$ of indecomposable objects in $\mathcal{C}$.
As an application, we describe Auslander-Reiten triangles, tilting
objects and exchange relations, and show a geometric method to construct the
Auslander-Reiten quiver of any cluster tilted algebra of type $D_n$,
using a result of \cite{BMR1}, in section \ref{sect 5}.

The author thanks Philippe Caldero, Frederic Chapoton and Dylan
Thurston for interesting discussions on the subject.
\end{section} 

\begin{section}{Preliminaries}\label{sect 1}

\begin{subsection}{Notation}\label{sect notation}
Let $k$ be an algebraically closed field. 
If $Q$ is a quiver, we denote by $Q_0$ the set of vertices and by $Q_1$
the set of arrows of $Q$. The {\em path algebra} of $Q$ over $k$ will be
denoted by $kQ$. It is of {\em finite representation type} if there is
only a finite  number of isoclasses of indecomposable modules. By
Gabriel's theorem  $kQ$ is of finite representation type if and only if
$Q$ is a Dynkin quiver, that is, the underlying graph of $Q$ is a
Dynkin diagram of type $A_n,D_n$ or $E_n$ \cite{G1}.

If $A$ is an algebra, we denote by $\textup{mod}\, A$ the category of
finitely generated right $A$-modules and by $\ind A$ a full
subcategory whose objects are a full set of representatives of the
isoclasses of indecomposable $A$-modules. 
Let $\DA=\mathcal{D}^b(\textup{mod}\, A)$ denote the derived category
of bounded complexes of finitely generated $A$-modules.
For further facts about $\textup{mod}\, (A)$ and $\DA$ we refer the
reader to \cite{ASS,ARS,G2,R}.

If $\mathcal{A}$ is an additive $k$-category then its {\em additive hull}
$\oplus \mathcal{A}$ is defined as follows:
The objects of $\oplus \mathcal{A}$ are direct sums of objects in
$\mathcal{A}$, morphisms $\oplus_i X_i \to \oplus_j Y_j$ are given
componentwise by morphisms $X_i\to Y_j$ of $\mathcal{A}$ and the
composition of morphisms is given by matrix multiplication.
\end{subsection}

\begin{subsection}{Translation quivers}
Following \cite{Riedtmann}, we define a {\em stable translation
 quiver} $(\zG,\tau)$ to be a quiver $\zG=(\zG_0,\zG_1)$ without loops
 together with a bijection $\tau$ (the {\em translation}) such that 
the number of arrows from $y \to x$ is equal to the number of arrows
 from $\tau\,x\to y$ for any $x,y\in\zG_0$. 
Given a stable translation quiver $(\zG,\tau)$, a polarization of $\zG$ is
 a bijection $\sigma:\zG_1\to\zG_1$ such that 
$\sigma(\za):\tau\,x\to y$  for every arrow $\za:y\to x\in \zG_1$.
If $\zG $ has no  multiple arrows, then there is a unique polarization.

Given a quiver $Q$ one can construct a stable translation quiver
$\mathbb{Z}Q$ as follows: $(\mathbb{Z}Q)_0=\mathbb{Z}\times Q_0$ and
the number of arrows in $Q$ from $(i,x)$ to $ (j,y)$ equals the number of
arrows  in $Q$ from $x$ to $y$ if $i=j$, and equals the number of arrows in $Q$
from
$y$ to $x$ if $j=i+1$, and there are no arrows otherwise. The
translation $\tau$ is defined by $\tau((i,x))=(i-1,x)$.

The {\em path category } of $(\zG,\tau)$ is the category whose  objects are
the vertices of $\zG$, and given $x,y\in\zG_0$, the $k$-space of
morphisms from $x$ to $y$ is given by the $k$-vector space with basis
the set of all paths from $x $ to $y$. The composition of morphisms is
induced from the usual composition of  paths.

The {\em mesh ideal} in the path category of $\zG$ is the ideal
generated by the {\em mesh relations}
\begin{equation}\nonumber
m_x =\sum_{\za:y\to x} \sigma(\za) \za.
\end{equation}  
The {\em  mesh category } $\mesh (\zG,\tau)$ of $(\zG,\tau)$ is the
quotient of the path 
category of $(\zG,\tau)$ by the mesh ideal.

Important examples of translation quivers are the Auslander-Reiten
quivers of the derived categories of hereditary algebras of finite
representation type. We shall need the following proposition.

\begin{prop}\label{prop 1}
Let $Q$ be a Dynkin quiver.
 Then 
\begin{enumerate}
\item for any quiver $Q'$ of the same Dynkin type as $Q$, the derived
 categories $\mathcal{D}^b kQ$ and $ 
 \mathcal{D}^b kQ' $ are  equivalent.
\item the Auslander-Reiten quiver of $\mathcal{D}^b kQ$
  is $\mathbb{Z}Q$.
\item the category $\ind\mathcal{D}^b kQ$ is equivalent
  to the mesh category of  $\mathbb{Z}Q$.
\end{enumerate}   
\end{prop}  

\begin{pf} See \cite[I.5]{H}. 
\qed
\end{pf}  
\end{subsection}

\begin{subsection}{Cluster categories }\label{sect cluster categories}
Let $\CD$ be the cluster category of type $D_n$, see \cite{BMRRT}. 
By definition, $\CD$
is the quotient of the derived category $\DA$ of a hereditary algebra
$A$ of type $D_n$ by the endofunctor
$F=\tau_{\DA}^{-1}\,[1]$, where
$\tau_{\DA}$ is the Auslander-Reiten
translation in $\DA$ and $[1]$ is the shift.
Thus the objects $\tilde M$ of $\CD$ are the orbits $\tilde M= 
(F^i\, M)_{i\in\ZZ}$ of
objects $M\in\DA$ and $\Hom_{\CD}(\tilde M,\tilde
N)=\oplus_{i\in\ZZ}\,\Hom_{\DA}(M,F^i\,N) $.
Let us denote the vertices of the Dynkin diagram of $D_n$ as follows
\[ 
\xymatrix@R0pt@C14pt{
&&&&& n \ar@{-}[dl]\\
1
&2\ar@{-}[l]\ar@{-}[r]
&\ldots \ar@{-}[r]
&(n-3)\ar@{-}[r]
&(n-2)\\
&&&&&(n-1)\ar@{-}[ul]
}
\]
and, for convenience, let us choose the algebra $A$ to be the path
algebra $kQ$ of the quiver 
\[\xymatrix@R0pt@C14pt{
&&&&&& n \\
Q=&1\ar[r]
&2\ar[r]
&\ldots \ar[r]
&(n-3)\ar[r]
&(n-2)\ar[ru]\ar[rd]\\
&&&&&&(n-1).
}
\]
Since, by Proposition \ref{prop 1}, the
Auslander-Reiten quiver of $\DA$ is the stable translation 
quiver $\ZZ Q$, the labels of $Q_0$ induce labels 
on the vertices of $\ZZ Q$ as usual (see Figure \ref{fig ar}):
\[ (\ZZ Q)_0=\{(i,j)\mid i\in \ZZ, \ j\in Q_0\}=\ZZ\times\{1,\ldots,n\}.\]
Moreover, we can identify the indecomposable objects of $\DA$ with the
vertices of $\ZZ Q$. Let $P_1$ be the indecomposable projective module
corresponding to the vertex $1\in Q_0$. Then, by defining the position
of $P_1$ to be $(1,1)$, we have a bijection
\[
 \pos : \ind \DA \to \ZZ\times\{1,\ldots,n\}.
\] 
In other terms, for $M\in\ind \DA$, we have $\pos (M)=(i,j)$ if and
only if $M=\tau_{\DA}^{-i}\,P_j$, where $P_j$ is the indecomposable
projective $A$-module at vertex $j$. 
The integer $j\in\{1,\ldots,n\}$ is called the {\em level} of $M$ and
will be denoted by $\level(M)$.

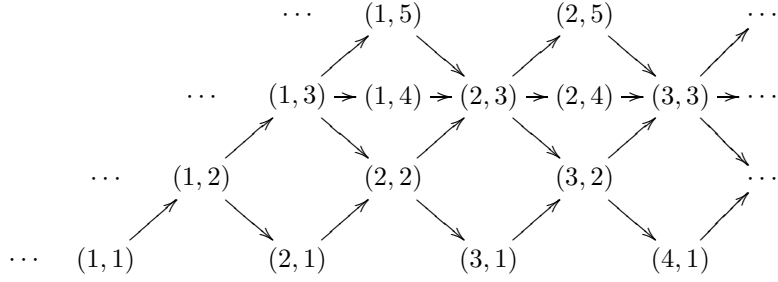
\begin{figure}
\[\xymatrix@R=15pt@C=8pt{
&&&\cdots&  (1,5)\ar[rd]&&{  (2,5)}\ar[rd]&&\cdots\\ 
&&\cdots&{(1,3)}\ar[ru]\ar[rd]\ar[r]
&{(1,4)}\ar[r]&{(2,3)}\ar[ru]\ar[rd]\ar[r]&
{(2,4)}\ar[r]&{(3,3)}\ar[ru]\ar[rd]\ar[r]&{\cdots}\\ 
&\cdots&{(1,2)}\ar[ru]\ar[rd]&&{(2,2)}\ar[ru]\ar[rd]
&&(3,2)\ar[ru]\ar[rd]&& \cdots\\  
\cdots&{{(1,1)}}\ar[ru]&&{{(2,1)}}\ar[ru]
&&{{(3,1)}}\ar[ru]&&{{(4,1)}}\ar[ru]
}\]
\caption{Labels of the Auslander-Reiten quiver of $\DA$ if $n=5$}\label{fig ar}
\end{figure}

If $\pos(M)=(i,j)$ with $j\in \{n-1,n\}$ then
let $M^-$ be the indecomposable object such that $\pos(M^-)=(i,j')$,
where $j'$ is the unique element in $\{n-1,n\}\setminus\{j\}$.
The structure of the module category of $A$ and of its derived
category is well known.
  In particular, we have the following result. 

\begin{lem}\label{lem 4} Let  $M\in\ind\DA$.
\begin{enumerate}
\item If $n$ is even, then $M [1]=\tau_{\DA}^{-n+1}\,M$.
\item If $n$ is odd, then 
\[M[1]=\left\{\begin{array}{ll} 
\tau_{\DA}^{-n+1}\,M &\textup{if  $\level(M)\le n-2$}\\
\\
\tau_{\DA}^{-n+1}\,M^- &\textup{if $\level(M)\in\{n-1,n\}$}
\end{array}\right. \]
\end{enumerate}   
\end{lem}    

\begin{pf}
It suffices to prove the statement in the case where $M=P$ is an
indecomposable projective $A$-module. Let $\nu_A$ denote the Nakayama
functor of $\textup{mod}\, A$. Then
$P[1]=\tau^{-1}_{\DA}\,\nu_A\,P$. Now the statement follows from
\cite[Proposition 6.5]{G2}.
\qed
\end{pf}  

We also define the position of indecomposable objects of the cluster
category $\CA=\DA/F$. The set
$\textup{mod}\, A \,\sqcup \, A[1]$ is a fundamental domain for the
cluster category (here $A[1]$ denotes the first shift of all
indecomposable projective $A$-modules).
Then, for $\tilde X\in\ind\CA$, we define $\pos (\tilde X)$ to be
$\pos(X)$ where $X$ is the unique element in the fundamental domain
such that $\tilde X$ is the orbit of $X$.

We recall now a well-known result, that uses the position
 coordinates to describe  
 the dimension of the space of morphisms between indecomposable
 $\CA$-objects. 
 Let $M,N\in\ind\CA$ such that $\pos(M)=(1,m)$ and $\pos(N)=(i,j)$.
 The next proposition follows
from the structure of the mesh category.

\begin{prop}\label{prop 2.2}
The dimension of $\Hom_{\CA}(M,N)$ is $0,1$ or $2$ and
it  can be characterized as follows:
\begin{enumerate}
\item If $m\le n-2$ then 
$ \dim\Hom_{\CA}(M,N)\ne 0 $ if and only if
 \[\left\{
\begin{array}{lccc}
&1\le i\le m &\& & i+j\ge m+1\\
\textup{or }&  m+1\le i \le n-1 &\& &  n\le i+j\le n+m-1
\end{array}  
\right.
\]
\item If $m\in\{n-1,n\}$ define $m'$ to be the unique element in
  $\{n-1,n\}\setminus\{m\}$.  Then 
$\dim\Hom_{\CA}(M,N)\ne 0$ if and only if 
\[\left\{
\begin{array}{lccccc}
& 2\le i\le n-1 &\&&
 i+j\ge n &\&&  j\le n-2\\
\textup{or } &1\le i\le n-1 &\&& j=m &\&& i \textup{ is odd}\\
\textup{or } &1\le i\le n-1 &\&& j=m' &\&& i \textup{ is even}
\end{array}  
\right.
\]
\end{enumerate}   
Moreover $\dim\Hom_{\CA}(M,N)=2$ if and only if
\[ 
\begin{array}{cccccccc} 
 2\le m\le n-2 &\&& 2\le i\le m &\&& 2\le j\le n-2 &\&& i+j\ge n.
\end{array}  
\]
\end{prop}

To illustrate this statement, we give an example where $n=6$ and $M$
is on position $(1,3)$. The following picture shows the
Auslander-Reiten quiver of $\CA$ (without arrows)  and where the label
on any vertex $N$ is $\dim\Hom_{\CA}(M,N)$. If this number is zero, we
write a dot instead of $0$. The two underlined
vertices should be identified, they indicate the position of
$M$. 
\[\begin{array}{cccccccccccccccccccccccccccc}
&.&&1&&1&&1&&.&&.&&.&&1&&1&&1&&.&\\
.&.&1&1&2&1&2&1&1&.&.&.&.&.&1&1&2&1&2&1&1&.&.\\
&\underline{1}&&1&&2&&1&&1&&.&&\underline{1}&&1&&2&&1&&1\\
.&&1&&1&&1&&1&&.&&.&&1&&1&&1&&1&&.\\
&.&&1&&.&&1&&.&&.&&.&&1&&.&&1&&.\\
\end{array}  
\]
\end{subsection}
\end{section} 


\begin{section}{The category of tagged edges}\label{sect 2}
\begin{subsection}{Tagged edges} \label{edges} 
Consider a regular polygon with $n\ge 3$ vertices and one puncture in its
center. 
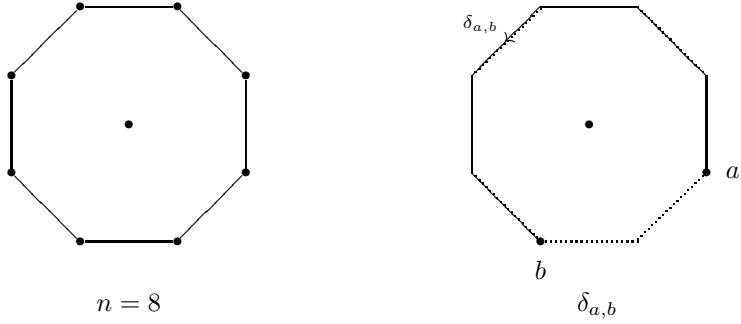
\[
\begin{array}{ccc}
 \xy/r4pc/:{\xypolygon8"A"{~<{}~>{-}{\scriptstyle\bullet}}}
*+{{\scriptstyle \bullet}},
\endxy 
&\hspace{2cm}&
\xy/r4pc/:{\xypolygon8"A"{~<{}~>{.}{}}}
*+{{\scriptstyle \bullet}},
\POS"A6"\drop{\scriptstyle\bullet}
\POS"A6"\drop{\begin{array}{c}\\ \\ b\end{array}}
\POS"A8"\drop{\scriptstyle\bullet}
\POS"A8"\drop{\qquad a} 
\POS"A8" \ar@{-}   "A1"
\POS"A1" \ar@{-}   "A2"
\POS"A2" \ar@{-}   "A3"
\POS"A3" \ar|-{\SelectTips{cm}{}\object@{>}}@{-}_{\zd_{a,b}}   "A4"
\POS"A4" \ar@{-}   "A5"
\POS"A5" \ar@{-}   "A6"
\endxy \\
n=8 &&\zd_{a,b}
\end{array}  \]


If $a\ne b$ are any two vertices on the boundary then let 
$\zd_{a,b}$
denote a path along the boundary from $a$ to $b$ in counterclockwise
direction which does not run through the same point twice.
Let $\zd_{a,a}$ denote a path along the boundary from $a$ to $a$ in
counterclockwise direction which goes around the polygon exactly
once, and such that $a$ is the only point through which $\zd_{a,a}$
runs twice.
For $a\ne b$, let $ | \zd_{a,b}| $ be the number of vertices on
the path $\zd_{a,b}$ (including $a$ and $b$), and let $|
\zd_{a,a}|  = n+1$.
  Two vertices $a,b$ are called {\em neighbors} if $| \zd_{a,b}| 
=2$ or $| \zd_{a,b}|   = n$, and
$b$ is the {\em counterclockwise neighbor} of $a$ if $|
\zd_{a,b}|  =2$.

An {\em edge} is a triple $(a,\za,b)$ where $a$ and $b$ are vertices
of the punctured polygon and $\za$ is a path from $a$ to $b$ such that 
\begin{itemize}
\item[(E1)] $\za$ is homotopic to $\zd_{a,b}$,
\item[(E2)] except for its starting point $a$ and its endpoint $b$, the path
  $\za$ lies in the interior of the punctured polygon,
\item[(E3)] $\za$ does not cross itself, that is, there is no point in the
  interior of the punctured polygon through which $\za$ runs twice,
\item[(E4)]  $| \zd_{a,b}|   \ge 3$.
\end{itemize}   
Note that condition (E1) implies that  if $a=b$ then $\za$ is not
homotopic to the constant path, and condition (E4) means that $b$ is
not the counterclockwise neighbor of $a$.

Two edges $(a,\za,b), (c,\zb,d)$ are {\em equivalent} if 
$a=c,\  b=d$  and $\za$ is homotopic to $ \zb$.
Let $E$ be the set of equivalence classes of edges.  Since the
homotopy class of an edge is already fixed by condition (E1), an
element of $E$ is uniquely determined by an ordered pair of vertices
$(a,b)$. We will therefore use the notation $M_{a,b}$ for the
equivalence class of edges $(a,\za,b)$ in $E$.

Define the set of {\em tagged edges} $\edges $ as follows:
\begin{equation}\nonumber
\edges=\{M_{a,b}^\ze \mid M_{a,b}\in E, \ze=\pm 1 \textup{ and
  $\ze=1$ if $a\ne b$}
\}
\end{equation} 
If $a\ne b$, we will often drop the exponent and write 
$M_{a,b}$ instead of $M_{a,b}^1$.
In other terms, for any ordered pair $(a,b)$ of vertices with $a\ne b$
and $b$ not the counterclockwise neighbor of $a$, there is 
exactly one tagged edge $M_{a,b}\in\edges$,
and, for any vertex
$a$, there are exactly two tagged edges $M_{a,a}^{-1}$ and $M_{a,a}^{1}$.
A simple count shows that there are $\frac{n!}{(n-2)!}-n+2n  =n^2$
elements in $\edges$. These tagged edges will correspond to the
indecomposable objects in the cluster category.

\begin{rem} In  pictures, the  tagged edges $M_{a,a}^\ze$ that start
  and end at the same vertex $a$ will not be represented  as loops
  around the puncture but as  lines from the vertex
  $a$ to the puncture. If $\ze=-1$ the line will have a tag on it and
  if $\ze=1$ there will be no tag. The reason for this is the
  definition of the crossing number in the next section.
\end{rem} 
\end{subsection}

\begin{subsection}{Crossing Numbers}\label{sect 1.2}

Let $M=M_{a,b}^\ze$ and $N=M_{c,d}^{\ze'}$ be in $\edges$ and let $\Delta^0$
  denote the interior of the punctured polygon. The {\em
  crossing number} $e(M,N)$ of $M$ and $N$ is the minimal number of
  intersections of $M_{a,b}$ and $M_{c,d}$ in $\Delta^0$. More precisely,
\begin{enumerate}
\item If $a\ne b$ and $c\ne d$, then 
\[e(M,N)=\min\{
  \textup{Card}(\za\cap\zb\cap \Delta^0 )\mid (a,\za,b)\in M_{a,b},\
  (c,\zb,d)\in M_{c,d} \}.\]  
\item If $a=b$ and $c\ne d$,  let $\za$ be  the
  straight line from $a$ to the puncture. Then
\[e(M,N)=\min\{ 
  \textup{Card}(\za\cap\zb\cap \Delta^0 )\mid (c,\zb,d)\in M_{c,d} \}.\] 
\item If $a\ne b$ and $c= d$,  let $\zb$ be  the
  straight line from $c$ to the puncture. Then
\[e(M,N)=\min\{ 
  \textup{Card}(\za\cap\zb\cap \Delta^0 )\mid (a,\za,b)\in M_{a,b} \}.\] 
\item If $a=b$ and $c=d$, that is $M=M_{a,a}^{\ze}$ and
  $N=M_{c,c}^{\ze'}$, then  
\[e(M,N)=\left\{
\begin{array}{ll}
1&\textup{if $a\ne c$ and $\ze\ne \ze'$}\\
0&\textup{otherwise}
\end{array}\right. \]
\end{enumerate}   
We say that $M$ {\em crosses } $N$ if $e(M,N)\ne 0$.
The next lemma follows immediately from the construction.

\begin{lem}\label{lem e0}
For any $M,N\in \ind \mathcal{C}$ we have $e(M,N)\in\{0,1,2\}$
\end{lem}

\end{subsection} 

\begin{subsection}{Elementary moves}\label{sect elementary moves}
Generalizing a concept from \cite{CCS}, we will now define elementary
  moves, which will correspond to irreducible morphisms in the cluster
  category. 

An elementary move  sends a tagged edge  $M_{a,b}^\ze\in \edges$ to
  another tagged edge
  $M_{a',b'}^{\ze'}\in\edges$ (i.e. it is an ordered pair $
  (M_{a,b}^{\ze},M_{a',b'}^{\ze'}) $ of tagged edges) satisfying certain
  conditions which we will define in four separate cases according to
  the relative position of $a$ and $b$.
Let $c$ (respectively $d$) be the counterclockwise neighbor of $a$
  (respectively $b$).
\begin{enumerate}
\item If $| \zd_{a,b}|  =3$, then there is precisely one elementary
  move    $M_{a,b}\mapsto M_{a,d}$.
\item If $4\le | \zd_{a,b}| \le n-1$, then 
there are precisely two elementary
  moves  $M_{a,b}\mapsto M_{c,b}$ and $M_{a,b}\mapsto M_{a,d}$.
\item If $| \zd_{a,b}|  =n$, then $d=a$ and
there are precisely three elementary
  moves  $M_{a,b}\mapsto M_{c,b}$, $M_{a,b}\mapsto M_{a,a}^1$ and
  $M_{a,b}\mapsto  M_{a,a}^{-1}$.
\item If $| \zd_{a,b}|    =n+1$, then $a=b$ and there is precisely one elementary move
  $M_{a,a}^\ze\mapsto M_{c,a}$.
\end{enumerate}  
Some examples of elementary moves are shown in Figure \ref{fig
  elementary moves}. 
\begin{figure}
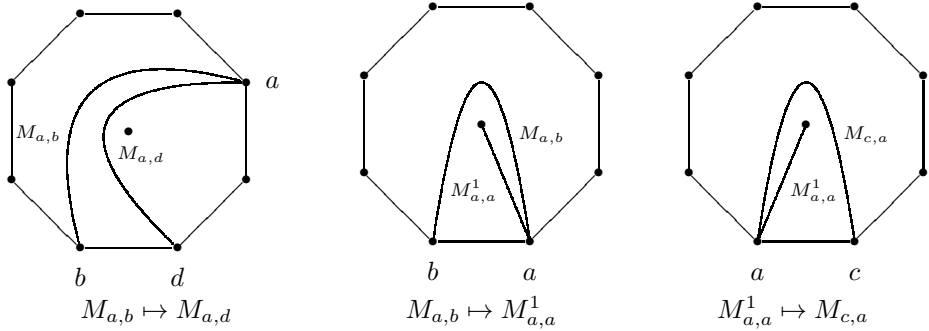
 
\[\begin{array}{ccc}
\xy/r4pc/:{\xypolygon8"A"{~<{}~>{-}{\scriptstyle \bullet}}}
*+{\scriptstyle\bullet},
\POS"A1"\drop{\qquad a}
\POS"A7"\drop{\begin{array}{c}\\ \\d \end{array}  } 
\POS"A6"\drop{\begin{array}{c}\\ \\b \end{array}  } 
\POS"A6" \ar@{}^(.3){M_{a,b}}@/^9ex/ "A1"
\POS"A7" \ar@{}_(.3){M_{a,d}}@/^10ex/ "A1",
\endxy
\quad 
&\quad
\xy/r4pc/:{\xypolygon8"A"{~<{}~>{-}{\scriptstyle \bullet}}}
*+{\scriptstyle\bullet},
\POS"A7"\drop{\begin{array}{c}\\ \\a \end{array}  } 
\POS"A6"\drop{\begin{array}{c}\\ \\b \end{array}  } 
\POS"A6" \ar@{}^(.8){M_{a,b}}@/^14ex/ "A7"
\POS"A7" \ar@{}^(.5){M_{a,a}^1} @/^0ex/ "A0",
\endxy 
\quad
&\quad
\xy/r4pc/:{\xypolygon8"A"{~<{}~>{-}{\scriptstyle \bullet}}}
*+{\scriptstyle\bullet},
\POS"A6"\drop{\begin{array}{c}\\ \\ a \end{array}  } 
\POS"A7"\drop{\begin{array}{c}\\ \\ c \end{array}  } 
\POS"A6" \ar@{}^(.8){M_{c,a}}@/^14ex/ "A7"
\POS"A6" \ar@{}_(.5){M_{a,a}^1} @/^0ex/ "A0",
\endxy \\
M_{a,b}\mapsto M_{a,d}
&M_{a,b}\mapsto M_{a,a}^1
&M_{a,a}^1\mapsto M_{c,a}
\end{array}\]
\caption{Examples of elementary moves}\label{fig elementary moves}
\end{figure}

\end{subsection}

\begin{subsection}{A triangulation of type $D_n$}\label{triangulation}
A {\em triangulation} of the punctured polygon is a maximal set of
 non-crossing  tagged edges. 
\begin{lem}\label{lem 3}
Any triangulation of the punctured polygon has $n$ elements.
\end{lem}  
\begin{pf} We prove the lemma by induction on $n$. If $n=3$, there
  exist,   up to rotation, symmetry and changing tags, three different
  triangulations (each having 3 tagged edges): 

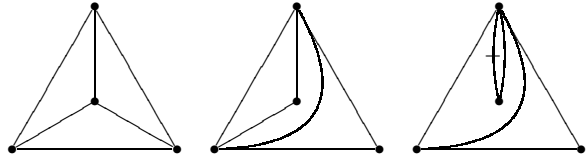
\[\begin{array}{ccc}
\xy/r3pc/:{\xypolygon3"A"{~<{}~>{-}{\scriptstyle \bullet}}}
*+{\scriptstyle\bullet},
\POS"A1" \ar@{-}   "A0"
\POS"A2" \ar@{-}   "A0"
\POS"A3" \ar@{-}   "A0"
\endxy
&\xy/r3pc/:{\xypolygon3"A"{~<{}~>{-}{\scriptstyle \bullet}}}
*+{\scriptstyle\bullet},
\POS"A1" \ar@{-}   "A0"
\POS"A2" \ar@{-}   "A0"
\POS"A2" \ar@{} @/_6ex/ "A1"
\endxy
&\xy/r3pc/:{\xypolygon3"A"{~<{}~>{-}{\scriptstyle \bullet}}}
*+{\scriptstyle\bullet},
\POS"A1" \ar@{-}@/^0.5ex/   "A0"
\POS"A1" \ar|-(0.5){\SelectTips{cm}{}\object@{+}}@{}@/_0.5ex/   "A0"
\POS"A2" \ar@{} @/_6ex/ "A1"
\endxy
\end{array}
\] 
Suppose now that $n>3$ and
  let $T$ be a triangulation. If $T$ contains an edge of the
  form $M=M_{a,b}$ with $a\ne b$ then this edge cuts the punctured
  polygon into two parts. Let $p$ be the number of vertices of the
  punctured polygon  that lie in
  the part that contains the puncture. Then $n-p+2$ vertices lie in
  the other part, because $a$ and $b$ lie in both. The set
  $T\setminus\{ M\}$ defines triangulations on both parts. If $p\ge 3$ then, by
  induction, the number of edges in $T\setminus \{M\}$ that lie in the part
  containing the puncture is equal to $p$. If $p=2$, then this number
  is also equal to $p$ (in this case the $2$ edges are either
  $(M_{a,a}^{\ze},M_{a,a}^{-\ze})$ or $(M_{a,a}^{\ze},M_{c,c}^{\ze})$
  or $(M_{c,c}^{\ze},M_{c,c}^{-\ze})$). Thus, in all cases the number
  of edges in $T\setminus \{M\}$ that lie in the part 
  containing the puncture is equal to $p$. 
On the other hand,  the
  number of edges in $T\setminus \{M\}$  that lie in other part equals the
  number of edges in a triangulation of an  $n-p+2$-polygon, and this
  number is equal to $n-p-1$. So the cardinality of  $T\setminus \{M\}$ is
  $n-1$, hence $T$ has $n$ elements.

  If $T$ does not contain an edge of the form $M_{a,b}$ with $a\ne
  b$, then $T$ contains exactly one element $M_{a,a}^{\ze}$ for each
  vertex $a$ and all edges have the same tag $\ze =1$ or $\ze=-1$.
  Clearly, $T$ has $n$ elements.
\qed
\end{pf}  

We construct now a particular triangulation  of the punctured
polygon.
 Let $a_1$ be a vertex on the boundary of the punctured polygon. Denote by
$a_2$ the counterclockwise neighbor 
 $a_1$ and, recursively, let $a_k$ be the
counterclockwise neighbor of $a_{k-1}$,
for any $k$ such that $2\le k \le n$.
Then the triangulation $T=T(a_1) $ is the  set of tagged edges
\[T =\{M_{a_1,a_1}^{1},M_{a_1,a_1}^{-1} \} \cup \{
    M_{a_1,a_{k}}\mid 3\le  k\le n \}.
\]
Here is an example for $T  $ in the case $n=8$.

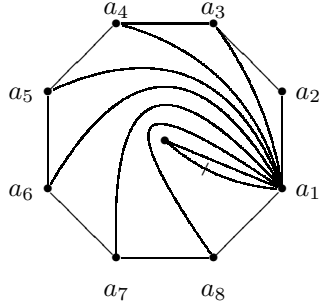
\[\xy/r4pc/: {\xypolygon8"A"{~<{}~>{-}{\scriptstyle\bullet}}},
*+{\scriptstyle\bullet},
\POS"A8"\drop{\begin{array}{c}\qquad a_1 \end{array}  } 
\POS"A7"\drop{\begin{array}{c}\\ \\ a_8 \end{array}  }
\POS"A3"\drop{\begin{array}{c} a_4 \\ \\ \end{array}  } 
\POS"A4"\drop{\begin{array}{c} a_5 \qquad  \end{array}  }
\POS"A1"\drop{\begin{array}{c}\qquad a_2  \end{array}  }  
\POS"A2"\drop{\begin{array}{c} a_3\\ \\  \end{array}  }  
\POS"A6"\drop{\begin{array}{c}\\ \\  a_7  \end{array}  }  
\POS"A5"\drop{\begin{array}{c} a_6 \qquad \end{array}  }  
\POS"A8" \ar@{-} @/^0ex/  "A0"
\POS"A8" \ar|-(0.6){\SelectTips{cm}{}\object@{|}}@{-}@/^1ex/ "A0"
\POS"A8" \ar@{-} @/_12ex/ "A7",
\POS"A8" \ar@{-} @/_11ex/ "A6",
\POS"A8" \ar@{-} @/_9ex/ "A5",
\POS"A8" \ar@{-} @/_6ex/ "A4",
\POS"A8" \ar@{-} @/_3ex/ "A3",
\POS"A8" \ar@{-} @/_1ex/ "A2"
\endxy 
\]

From the construction it is clear that the  elementary moves between
the elements of $T $ are precisely the elementary moves  
$M_{a_1,a_k}\mapsto M_{a_1,a_{k+1}}$ for all $k$ with $3\le k\le n-1$ and
$M_{a_1,a_n}\mapsto M_{a_1,a_{1}}^\ze$ with $\ze = \pm 1$.
We can associate a quiver $Q_T $ to $T $ as follows: The
  vertices  of
 $Q_T$ are the elements of  $T $ and the arrows are the
  elementary moves between the elements of $T $.
By the above observation, $Q_T $ is the following
 quiver of type $D_n$ (for $n=3$ set $D_3=A_3$):
\[
\xymatrix@R0pt@C14pt{
&&&&M_{a_1,a_1}^{1}\\
M_{a_1,a_3} \ar[r] &M_{a_1,a_4} \ar[r] &\cdots  \ar[r] &
M_{a_1,a_n}\ar[rd]\ar[ru] \\
&&&&M_{a_1,a_1}^{-1}.
}
\]
Note that $Q_T$ is isomorphic to the quiver $Q$ that we have chosen to
construct the cluster category in section \ref{sect cluster categories}.
\end{subsection}

\begin{subsection}{Translation} \label{sect tau}
We define the {\em translation} $\tau$ to be the following  bijection
$\tau:\edges\to\edges$:
Let $M_{a,b}^\ze\in \edges$ and let $a',b'$ be the two vertices such
that  $a$ (respectively $b$) is the
  counterclockwise neighbor of $a'$ (respectively $b'$). 
\begin{enumerate}
\item If $a\ne b$ then $\tau  M_{a,b}=M_{a',b'}$.
\item If $a= b$ then  $\tau  M_{a,a}^\ze=M_{a',a'}^{-\ze}$, for
  $\ze=\pm 1$.
\end{enumerate} 
The next lemma follows immediately from the definition of $\tau$.
\begin{lem}\label{lem 1} 
\begin{enumerate}
\item If $n$ is even then $\tau^n =\textup{id}$.
\item  If $n$ is odd then $\tau^n \,M_{a,b}^\ze=\left\{
\begin{array}{ll} 
M_{a,b}^\ze &\textup{if $a\ne b$},\\
\\
M_{a,a}^{-\ze} &\textup{if $a= b$}.
\end{array}  \right.$
\end{enumerate}   
\end{lem}  

We will need the following Lemma when we define the category of tagged
edges.

\begin{lem}\label{lem 2}
Let $M_{a,b}^\zl,\ \ M_{c,d}^{\ze}\in\edges$. Then there is an elementary move
$M_{a,b}^{\zl} \mapsto M_{c,d}^{\ze} $ if and only if there is an elementary
move
$\tau  M_{c,d}^{\ze} \mapsto M_{a,b}^{\zl}$.
\end{lem}  

\begin{pf} Let $ M_{c',d'}^{\ze'} = \tau \,  M_{c,d}^{\ze}$.
Suppose that there is an elementary move
$M_{a,b}^{\zl} \mapsto M_{c,d}^{\ze} $.  
Then either $a=c$  or $b=d$.
\begin{enumerate}
\item If $a=c$, then $b=d'$, and $| \zd_{c,d}| = | \zd_{c',d'}|   =
  |\zd_{a,b}|  +1$.  In particular, $a\ne
  b$ and $ | \zd_{c',d'}|  \ge 4$. 
Now, if $4\le  |  \zd_{c',d'}| \le n$, then $c'\ne d'$ and, by (\ref{sect
  elementary moves}.2, \ref{sect elementary moves}.3), there is an
  elementary move  $\tau\,M_{c,d}=M_{c',b}\mapsto M_{c,b}=M_{a,b}$.
On the other hand, if $ |  \zd_{c',d'}| =n+1$ then $c'=d'=b$ and
 $a$ is the
  counterclockwise neighbor of $b$ and thus, by  (\ref{sect
  elementary moves}.4), there is an elementary move
  $\tau\,M_{c,d}^\ze=M_{b,b}^{-\ze} \mapsto M_{a,b}$.
\item If $b=d$ then $a=c'$ and  $| \zd_{c,d}|=| \zd_{c',d'}|  = |
  \zd_{a,b}|  -1$. In particular, $ |   \zd_{c',d'}|  \le n$.
Now, if $3\le  |  \zd_{c',d'}| \le n-1$, then, by  (\ref{sect
  elementary moves}.1, \ref{sect elementary moves}.2), there is an
  elementary move  $\tau\,M_{c,d}=M_{a,d'}\mapsto M_{a,b}$.
On the other hand, if $ |  \zd_{c',d'}| =n$ then $a=b=d$ and $b$ is the
  counterclockwise neighbor of $d'$ and thus, by  (\ref{sect
  elementary moves}.3), there is an elementary move
  $\tau\,M_{c,d}^\ze=M_{a,d'} \mapsto M_{a,b}^\zl$, for $\zl=\pm 1$.
\end{enumerate}   
Thus there is an elementary move  $ \tau\,M_{c,d}^{\ze}\mapsto M_{a,b}^\zl$.
By symmetry, the other implication also holds.
\qed
\end{pf}

\end{subsection}

\begin{subsection}{The category $\mathcal{C}$}\label{sect C}

We will now define a $k$-linear additive category $\mathcal{C}$ as
follows. The objects are direct sums of tagged edges in $\edges$. By
additivity, it is enough to define morphisms between tagged edges. The
space of morphisms from a tagged edge $M$ to a tagged edge
$N $ is a quotient of the vector space over $k$ spanned
  by sequences of elementary moves starting at $M$
and ending at $N$. The subspace which defines the
  quotient is spanned by the so-called {\em mesh relations}.

For any tagged edge $X\in \edges$, we define the mesh relation
\[ m_X=\sum_{Y\stackrel{\za}{\mapsto} X} \tau
X\stackrel{\zs(\za)}{\mapsto} Y \stackrel{\za}{\mapsto} X\] 
where  the sum is over all elementary moves
$\za$ that send a tagged edge $Y$ to $X$,
and $\zs(\za)$ denotes the elementary move given by Lemma \ref{lem
  2}. 
In other terms, if there is only one $Y\in\edges$ such that there is
an elementary move  $Y\stackrel{\za}{\mapsto} X$ then the composition
of morphisms $\tau\,X
\stackrel{\zs(\za)}{\mapsto} Y\stackrel{\za}{\mapsto} X$  is zero, and
if there are several $Y_1,\ldots ,Y_p\in \edges$ with elementary moves
$Y_i \stackrel{\za_i}{\mapsto} X$ then the composition of morphisms $\tau\,X
\stackrel{\zs(\za_1)}{\mapsto}Y_1\stackrel{\za_1}{\mapsto}  X$
is equal to the sum of all compositions of morphisms 
$\tau\,X\stackrel{\zs(\za_i)}{\mapsto}Y_i\stackrel{\za_i}{\mapsto}  X$,
 $i=2,\ldots,n$.

More generally, a mesh relation is an equality between two sequences
of elementary moves which differ only in two consecutive elementary
moves by such a relation $m_X$.

We can now define the set of morphisms from a tagged edge $M$ to a
tagged edge $N$ to be the quotient of the vector space over $k$
spanned by sequences of elementary moves from $M$ to $N$ by the
subspace generated by mesh relations.
\end{subsection} 
\end{section} 


\begin{section}{Equivalence of categories}\label{sect 3}
In this section, we will prove the equivalence of the category $\mathcal{C}$
and the cluster category $\CA$.

\begin{subsection}{Translation quiver}\label{sect tq}

First we construct a stable translation quiver in our
situation. Let $\zG $ be the following quiver:
 The set of vertices  is the set $\mathbb{Z}\times \edges$
 where $\edges$ is the set of tagged edges of section \ref{edges}.
In order to define the set of arrows, let us fix a vertex
$a_1$ on the boundary of the punctured polygon and let  $T=T(a_1) $ be
the triangulation  of section
\ref{triangulation}. 
Given $M,N\in\edges$, there is an arrow $(i,M)\to
 (i,N)$ in $\zG_1$ if
there is an elementary move $M\mapsto N$ and either $N\notin T $
or $M$ and $N$ are both in $T $,  and 
 there is an arrow
$(i,M)\to (i+1,N)$ in $\zG_1$ if  
there is an elementary move $M\mapsto N$ and $N\in T $ and $M\notin T$.
Note that $\zG $ has no loops and no  multiple arrows.

The translation $\tau : \edges \to \edges $  defined in section
\ref{sect tau} induces a bijection $\tau_{_\zG}$ on $\zG_0$ by
\[ \tau_{_\zG} (i,M)=\left\{ 
\begin{array}{ll} 
(i,\tau M) &\textup{if $M\notin T $}\\
(i-1,\tau M) &\textup{if $M\in T $.}
\end{array}  \right.
\]
 Then it 
follows from Lemma \ref{lem 2} that $(\zG ,\tau_{_\zG}) $ is a stable
translation quiver.
Let $\mesh (\zG ,\tau_{_\zG})$ denote the mesh category of $(\zG ,\tau_{_\zG})$.

Now consider  the vertex subset  $0\times
T \subset \zG_0$. This set contains exactly one element of each
$\tau_{_\zG}$-orbit of $\zG $. Moreover,  
the full subquiver induced by  $0\times T $  is the quiver
$Q=Q_T $ of sections \ref{sect cluster categories} and
\ref{triangulation}, a quiver of 
type $D_n$, and $\zG $ is the translation quiver $\ZZ Q $. Let $A$ be the path
algebra of $Q$, as 
before, and let us use the notation $(\ )^-$ of section \ref{sect
  cluster categories}.  
Then, by Proposition \ref{prop 1}, we have the following
result.

\begin{prop}\label{prop 2}
There is an equivalence of
categories 
\[\varphi : \mesh(\zG ,\tau_{_\zG}) \stackrel{\sim}{\longrightarrow} \ind\DA\]
 such that  \begin{enumerate}
\item $\varphi$ maps the elements of $0\times T $
to the projective $A $-modules,
\item $\varphi \circ \tau_{_\zG} = \tau_{\DA} \circ \varphi$,
\item If $\,\level(\varphi (\ell,M_{a,b}^\ze))=j$ then
\[ j= | \zd_{a,b}|   -2 \textup{ if $a\ne b$ \quad and \quad }
j \in \{n-1,n\} \textup{ if $a=b$},
\]
\item If $M=(\ell,M_{a,a}^\ze)$ and $M^-=(\ell,M_{a,a}^{-\ze})$ then
\[(\varphi(M))^-=\varphi(M^-).\]
\end{enumerate}     
\end{prop}  

The equivalence $\varphi $ induces an equivalence of categories $\varphi':\oplus \mesh
(\zG ,\tau_{_\zG}) \stackrel{\sim}{\to} \DA$, where $\oplus \mesh
(\zG ,\tau_{_\zG})$ denotes the additive hull of $\mesh
(\zG ,\tau_{_\zG})$.
In particular, $\oplus \mesh
(\zG ,\tau_{_\zG})$ has a triangulated structure such that $\varphi'$ is
an equivalence of triangulated categories. 
Let $\rho$ be the endofunctor of $\oplus \mesh
(\zG ,\tau_{_\zG})$ given by the full counterclockwise rotation of the
punctured polygon. Thus, if $(i,M)\in \mesh
(\zG ,\tau_{_\zG})$ then $\rho(i,M)=(i+1,M)$.
The next result follows from Lemma \ref{lem 1}.

\begin{lem}\label{lem rho}
\begin{enumerate}
\item If $n $ is even then $\rho=\tau_{_\zG}^{-n}$.
\item If $n$ is odd and $M=(\ell,M_{a,b})$ with $a\ne b$ then
  $\rho(M)=\tau_{_\zG}^{-n}(M)$. 
 \item If $n$ is odd  and $M=(\ell,M_{a,a}^\ze)$ and
   $M^-=(\ell,M_{a,a}^{-\ze})$ 
then  $\rho(M)=\tau_{_\zG}^{-n}(M^-)$. 
\end{enumerate}   
\end{lem}

\end{subsection}

\begin{subsection}{Main result}
The category of tagged edges $\mathcal{C}$ of section \ref{sect C} is
 the quotient category of the triangulated
category $\oplus \mesh
(\zG ,\tau_{_\zG}) $ by the endofunctor $\rho$.
Hence the objects of $\mathcal{C}$ are the $\rho$-orbits $\tilde M
=(\rho^i M)_{i\in\ZZ}$ of objects $M$ in $\oplus \mesh
(\zG ,\tau_{_\zG})$, and  $\Hom_{\mathcal{C}}(\tilde M,\tilde N) = \oplus
_i\Hom(M,\rho^i M)$ where the Hom-spaces on the right are taken in the
category $\oplus \mesh (\zG ,\tau_{_\zG})$.
We are now able to show our main result.

\begin{thm}\label{thm 1}
There is an equivalence of categories 
\[\overline{\varphi}:\mathcal{C} \longrightarrow \CA\]
between the category of tagged edges $\mathcal {C}$ and the cluster
category $\CA$ which sends the $\rho$-orbit of the
triangulation $0\times T $ to  the $F$-orbit of  the projective
$A$-modules. 
\end{thm}

\begin{pf} 
Let $\varphi':\oplus\mesh(\zG,\tau_{_\zG})\to \DA$ be the 
equivalence of section \ref{sect tq}.
Since $C= \oplus\mesh(\zG,\tau_{_\zG})/\rho$ and
$\CA=\DA/F$, we only have to show that 
$\varphi(\rho(M))=F(\varphi(M))$ for any $M\in
  \mesh(\zG ,\tau_{_\zG})$.

If $n$ is even then $F\,\varphi=\tau_{\DA}^{-1} \,[1]\,\varphi =
\tau_{\DA}^{-n}\,\varphi$, by Lemma \ref{lem 4}.
On the other hand, $\varphi\,\rho=\varphi\,\tau_{_\zG}^{-n}$, by Lemma
\ref{lem rho}, and then the statement follows from Proposition
\ref{prop 2}.

Suppose now that $n$ is odd.
Let $M=(\ell,M_{a,b}^\ze)\in\mesh(\zG,\tau_{_\zG})$. 
If $a\ne b$ then $\level(\varphi(M))= | \zd_{a,b}| -2\le n-2$,
by Proposition \ref{prop 2}, and, therefore,  
$F(\varphi(M)) = \tau_{\DA}^{-n}\,\varphi (M)$, by  Lemma \ref{lem
  4}. 
On the other hand, Lemma \ref{lem rho} implies that 
$\varphi\,\rho(M)=\varphi\,\tau_{_\zG}^{-n}(M)$, and  
the statement follows from Proposition
\ref{prop 2}.
Finally, suppose that $a=b$, that is, $M=(\ell,M_{a,a}^\ze)$ and
$\level(\varphi(M))\in\{n-1,n\}$. Let $M^- = (\ell,M_{a,a}^{-\ze})$,
then, by Proposition \ref{prop 2}, we have 
$(\varphi(M))^-=\varphi(M^-)$ and thus
$F(\varphi(M)) = \tau_{\DA}^{-n}\,\varphi (M^-)$,
by  Lemma \ref{lem 4}. 
On the other hand, Lemma \ref{lem rho} implies that 
$\varphi\,\rho(M)=\varphi\,\tau_{_\zG}^{-n}(M^-)$, and  again
the statement follows from Proposition
\ref{prop 2}.
\qed
\end{pf}

\begin{rem}
It has been shown in \cite{K} that the cluster category $\CA$ is
triangulated. The shift functor $[1]$ of this triangulated structure
is induced by the shift functor of $\DA$ and $[1]=\tau_{\CA}$, by
construction.  

Thus the category of tagged edges is also triangulated and its shift
functor is equal to $\tau$. In particular, we can define the $\Ext^1$
of two objects $M,N\in\ind \mathcal{C}$ as
\[\Ext^1_{\mathcal{C}}(M,N)=\Hom_{\mathcal C}(M,\tau\,N)\]
We study $\Ext^1$ in the next section.
\end{rem} 

\end{subsection} 
\end{section}


\begin{section}{Dimension of $\Ext^1$}\label{sect 4}

We want to translate the statement of Proposition \ref{prop 2.2} 
in the category $\mathcal{C}$. For an
element $N\in\ind\mathcal{C}$ let $\pos(N)=(i,j)$ be the position of the
corresponding indecomposable object $\overline{ \varphi}(N)$ in
$\CA$ under the equivalence $\overline{ \varphi}$ of Theorem
\ref{thm 1}. 
For two vertices $a,b$ 
on the boundary of the punctured polygon, define the closed interval $[a,b]$
 to be the set of all vertices that lie
    on the counterclockwise path from $a$ to $b$ on the boundary. 
The open interval $]a,b[$ is $[a,b]\setminus\{a,b\}$. 
Recall that $e(M,N) $ denotes the crossing number of tagged edges $M,N$
(see section \ref{sect 1.2}).

\begin{lem}\label{lem e1} Let $M_{a,b}\in\ind\mathcal{C}$ with $a\ne
    b$ such that   $\pos(M_{a,b})=(1,m)$.
    Let $M_{x,y}^\ze\in\ind\mathcal{C}$ be arbitrary and denote $\pos(M_{x,y}^\ze)$ by
    $(i,j)$. Then 
\begin{enumerate}
\item 
$e(\tau M_{a,b},M_{x,y}^\ze)=1$ if and only if one of the
  following conditions hold: 
\begin{eqnarray*}
1\le i\le m
&\&& m+1 \le i+j \le  n-1  \label{1''}\\
m+1 \le i\le n-1 
&\&& n \le i+j \le  n+m-1 \label{2''}\\
1\le i\le m 
&\&& n-1\le j\le n.\label{3''}
\end{eqnarray*} 
\item 
$e(\tau M_{a,b},M_{x,y}^\ze)=2$ if and only if
\[\begin{array}{ccccccc}
m\ge 2 &\&& 2\le i\le m&\&& n\le i+j &\&& 2\le j\le n-2.
\end{array} \]
\end{enumerate}   
\end{lem}

\begin{pf} Let $M_{a,b}$ be as in the lemma. By Proposition \ref{prop
    2}, we have  
  $1\le m = | \zd_{a,b}|  -2\,\le n-2$. We 
 use the notation $\tau \,M_{a,b}=M_{a',b'}$ and $\tau^2
  \,M_{a,b}=M_{a'',b''}$. 
Let $M_{x,y}^\ze\in\ind\mathcal{C}$ be such that
  $\pos(M_{x,y}^\ze)=(i,j)$.
Then 
$e(\tau M_{a,b},M_{x,y}^\ze)=1$ if and only if one of the following
  conditions hold: 
\begin{eqnarray} 
 x\in [a,b''] &\&& y\in [b,a'] \label{1}\\
 x\in [b',a''] &\&& y\in [a,b''] \label{2}\\
 x=y\in [a,b'']\label{3}.
\end{eqnarray}    
Using the fact that, by Proposition \ref{prop 2}, we have  
$i= | \zd_{a,x}|  $, and if $x\ne y$ then 
$j=| \zd_{x,y}|  -2$  and $i+j =
| \zd_{a,y}|  -1$, we can rewrite these conditions in terms 
  of $i$ and $j$ as follows:
\begin{eqnarray} 
1\le i\le | \zd_{a,b''}|   
&\&& | \zd_{a,b}|  -1 \le i+j \le  | \zd_{a,a'}|  -1 \label{1'}\\
 | \zd_{a,b'}|  \le i\le | \zd_{a,a''}|   
&\&& | \zd_{a,a}|  -1 \le i+j \le  n
	 +| \zd_{a,b''}|  -1 \label{2'}\\
1\le i\le | \zd_{a,b''}|   
&\&& n-1\le j\le n,\label{3'}
\end{eqnarray} 
where  (\ref{1})$\ssi$(\ref{1'}),  (\ref{2})$\ssi$(\ref{2'}) and
(\ref{3})$\ssi$(\ref{3'}). 
Now the first statement of the lemma follows simply by counting the
vertices on the 
various boundary paths that appear in (\ref{1'})-(\ref{3'}).

On the other hand $e(\tau\,M_{a,b},M_{x,y}^\ze)=2$ if and only if 
\[ 
\begin{array}{ccccccc}
x\in \,]a,b'[ &\&& y\in\, ]a',x[ \,.
\end{array}   
\]
Note that $x$ cannot be equal to $a$ because otherwise $a=x=y$ and, then by
definition of the crossing number, $e(\tau\,M_{a,b}^\ze,M_{x,y}^\ze)\le 1$. In particular,
we have $m= | \zd_{a',b'}|  -2\ge 2$.
Again, using the fact that 
$i=| \zd_{a,x}| $,  and if $x\ne y$ then $j=| \zd_{x,y}|  -2$ and $i+j =
	     n+| \zd_{a,y}|  -1$ 
we see that $e(\tau\,M_{a,b},M_{x,y}^\ze)=2$ if and only if 
\[ 
\begin{array}{cccccccc}
m\ge 2&\&& 2\le i\le | \zd_{a,b'}|  -1 &\&& n+0\le i+j &\&&2\le j\le n-2.
\end{array}   
\]
Now, the second statement of the lemma  follows because $| \zd_{a,b'}|  -1=m$.
\qed
\end{pf}

\begin{lem}\label{lem e2}
 Let $M_{a,a}^\ze\in\ind\mathcal{C}$ such that  $\pos(M_{a,a}^\ze)=(1,m)$
    with $m\in\{n-1,n\}$. Define $m'$ to be the unique element in
  $\{n-1,n\}\setminus\{m\}$.
    Let $M_{x,y}^{\ze'}\in\ind\mathcal{C}$ be arbitrary and denote
    $\pos(M_{x,y}^{\ze'})$ by 
    $(i,j)$. Then 
\begin{enumerate}
\item 
$e(\tau M_{a,a}^\ze,M_{x,y}^{\ze'})=1$ if and only if one of the
  following conditions hold: 
\begin{eqnarray*}
 2\le i\le n-1 &\&&
 i+j\ge n \quad\&\quad  j\le n-2\\
 1\le i\le n-1 &\&& j=\left\{\begin{array}{ll}
m' &  \textup{if $i$ is odd}\\
m&  \textup{if $i$ is even.}
\end{array}\right.
\end{eqnarray*}  
\item  $e(\tau M_{a,a}^\ze,M_{x,y}^{\ze'})=0$ otherwise.
\end{enumerate}
\end{lem}  
\begin{pf} Let $M_{a,a}^\ze$ be as in the lemma and let $\tau \,M_{a,a}^\ze
  =M_{a',a'}^{-\ze}$ and  $\tau^2\, M_{a,a}^\ze = M_{a'',a''}^\ze$.
Let $M_{x,y}^{\ze'}\in\ind\mathcal{C}$ be such that
  $\pos(M_{x,y}^{\ze'})=(i,j)$.
By  definition of the crossing number, we have $e(\tau
M_{a,a}^\ze,M_{x,y}^{\ze'})=1$ if 
  and only if 
one of the following conditions hold:
\begin{eqnarray*}
x\in \,]a,a'']&\&& y\in [a,x[\label{a}\\
x=y\ne a' &\&& \ze\ne\ze' \label{b}
\end{eqnarray*}  
We can rewrite these conditions in terms 
  of $i$ and $j$ as follows: 
\begin{eqnarray*}
2\le i \le | \zd_{a,a''}|   &\&& n+0\le i+j  \quad\&\quad j\le
| \zd_{x,x}|  -3 \label{a'}\\
1\le i\le n-1  &\&& j=\left\{ 
\begin{array}{ll}
m'&\textup{ if $i$ is odd}\\
m &\textup{ if $i$ is even.}
\end{array}\right.  
\end{eqnarray*}  
Counting the vertices on the boundary paths yields the lemma.
\qed
\end{pf}  

\begin{thm}\label{thm 2}
Let $M,N\in\ind\mathcal{C}$. Then the dimension of
$\Ext_{\mathcal{C}}^1(M,N)$ is equal to the crossing number $e(M,N)$
of $M$ and $N$. 
\end{thm}  

\begin{pf} 
Combining Proposition \ref{prop 2.2} with the Lemmas \ref{lem e0},\
\ref{lem e1} and \ref{lem e2}, we see that 
the dimension of $\Hom_{\mathcal{C}}(M,N)$ is equal to $e(\tau\, M,N)$, for all
$M,N\in \ind \mathcal{C}$ (we may assume without loss of generality
that $\pos(M)=(1,m)$, for some $m\in\{1,\ldots,n\}$). Hence $e(M,N)=
\dim\Hom_{\mathcal{C}}(\tau^{-1}\,M,N)=\dim\Hom_{\mathcal{C}}(M,\tau\,N)$. 
On the other hand,  $\Ext_{\mathcal{C}}^1(M,N)=\Hom_{\mathcal{C}}(M,\tau N)$ by definition, and
the result follows.
\qed
\end{pf}  

\begin{cor}\label{cor 1}
Let $M,N\in\ind\mathcal{C}$. Then
\[
\dim\Ext_{\mathcal{C}}^1(M,N)=\dim\Ext_{\mathcal{C}}^1(N,M).\] 
\end{cor}  

\end{section}



\begin{section}{Applications}\label{sect 5}

\begin{subsection}{Auslander-Reiten-triangles}
By \cite{K} and Theorem \ref{thm 1}, 
$\mathcal{C}$ is a triangulated category with
Auslander-Reiten-triangles
\[\tau \,M \to L \to M \to \tau\, M [1]=\tau^2\,M .\]
We have described $\tau$ already, let us describe $L$ now.
\begin{enumerate}
\item Suppose that  $\tau\,M =M_{a,b}$ with $a\ne b$. Then
  $M=M_{c,d}$ where $c$ (respectively  $d$) is the counterclockwise
  neighbor of $a$ (respectively $b$). In particular, $c\ne b$ because
  $|\zd_{a,b}| \ge 3$.

If  $a\ne d$, then  $L$ is the direct sum of $M_{c,b}$
and $M_{a,d}$ and if $a=d$ then  $L$ is the direct sum
of $M_{c,b}$, $M_{a,a}^{1}$ and $M_{a,a}^{-1}$.

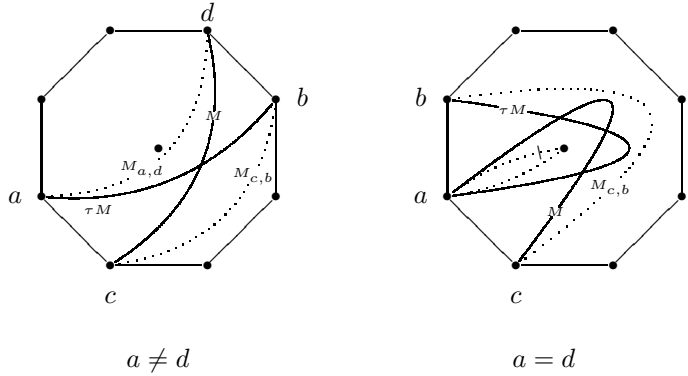
\[
\begin{array}{ccc}
 \xy/r4pc/: {\xypolygon8"A"{~<{}~>{-}{\scriptstyle\bullet}}},
*+{\scriptstyle\bullet},
\POS"A1"\drop{\begin{array}{l} \qquad b  \end{array}  }  
\POS"A2"\drop{\begin{array}{r} d\\ \\   \end{array}  }  
\POS"A6"\drop{\begin{array}{l}\\ \\  c  \end{array}  }  
\POS"A5"\drop{\begin{array}{r} a \qquad \end{array}  }  
\POS"A6" \ar@{-}|(0.7){\scriptscriptstyle  M}@/_{4ex}/  "A2",
\POS"A6" \ar@{.}|(0.7){\scriptscriptstyle M_{c,b}}@/_{4ex}/ "A1",
\POS"A5" \ar@{.}|(0.4){\scriptscriptstyle M_{a,d}\ }@/_{5ex}/ "A2",
\POS"A5" \ar@{-}_(0.2){\scriptscriptstyle \tau M}@/_{3ex}/ "A1"
\endxy 
&\quad&
 \xy/r4pc/: {\xypolygon8"A"{~<{}~>{-}{\scriptstyle\bullet}}},
*+{\scriptstyle\bullet},
\POS"A4"\drop{\begin{array}{l}b \qquad   \end{array}  }  
\POS"A6"\drop{\begin{array}{l}\\ \\  c  \end{array}  }  
\POS"A5"\drop{\begin{array}{r} a \qquad \end{array}  }  
\POS"A6" \ar@{-}|(0.1){\scriptscriptstyle M}@/_{16ex}/  "A5",
\POS"A5" \ar@{-}|(0.9){\scriptscriptstyle \tau M}@/_{16ex}/ "A4",
\POS"A5" \ar@{.}@/_{0.75ex}/"A0",
\POS"A5" \ar|-(0.8){\SelectTips{cm}{}\object@{|}}@{.}@/^{0.75ex}/"A0", 
\POS"A6" \ar@{.}|(0.2){\scriptscriptstyle M_{c,b}}@/_{16ex}/ "A4"
\endxy \\
\\
a\ne d &&a=d
\end{array} 
\]
\item Suppose that $\tau\, M=M_{a,a}^{-\ze}$. Then $M=M_{c,c}^{\ze}$,
  where $c$  is the counterclockwise
  neighbor of $a$. Then $L$ is indecomposable, $L=M_{c,a}$.

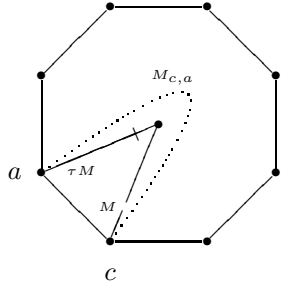
\[\xy/r4pc/: {\xypolygon8"A"{~<{}~>{-}{\scriptstyle \bullet}}},
*+{\scriptstyle\bullet},
\POS"A6"\drop{\begin{array}{l}\\ \\  c  \end{array}  }  
\POS"A5"\drop{\begin{array}{r} a \qquad  \end{array}  }  
\POS"A5" \ar@{.}@/^14ex/^(0.4){\scriptscriptstyle M_{c,a}}  "A6",
\POS"A6" \ar@{-}|(0.3){\scriptscriptstyle M\hspace{3ex}} "A0",
\POS"A5"\ar|-(0.8){\SelectTips{cm}{}\object@{+}}
@{-}|(0.05){\hspace{6ex}\scriptscriptstyle  \tau M}  "A0",
\endxy 
\]
\end{enumerate}   
Note that there are irreducible morphisms $\tau M \to L$ and $L\to M$,
and $\tau M \to L\to M$ is a mesh.
\end{subsection}

\begin{subsection}{Tilting objects and exchange pairs}
A tilting object $T$ in the category $\mathcal{C}$ is a maximal set
of non-crossing tagged edges, that is, $T$ is a triangulation of the
punctured polygon.
Let $M,N\in\edges$. 
If $e(M,N)=1$ then, by
\cite[Theorem 7.5]{BMRRT}, $M,N$ form an {\em exchange pair}, that is,
there exist two triangulations $T$ and $T'$ 
such that $M\in T$, $N\in T'$ and
$T\setminus\{M\}=T'\setminus\{N\}$. The edges $M$ and $N$ are the
``diagonals'' in a generalized quadrilateral in $T\setminus\{M\}$, see
Figure \ref{fig exchange}.
 The triangulation $T'$ is obtained
from the triangulation $T$ by ``flipping'' the diagonal $M$ to the
diagonal $N$.

Let $x_M,x_N$ be the corresponding cluster variables in the cluster
algebra, then the exchange relation is given by
\[x_M\,x_N=\prod_i x_{L_i} +\prod_{i} x_{L_i'}\]
where the products correspond to ``opposite'' sides in the generalized
quadrilateral and can have one, two or three factors, see Figure
\ref{fig exchange}. 
\begin{figure}
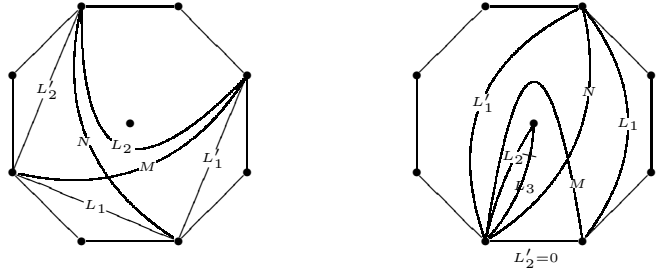

\[
\begin{array}{ccc}
\xy/r4pc/:{\xypolygon8"A"{~<{}~>{-}{\scriptstyle\bullet}}}
*+{{\scriptstyle \bullet}},
\POS"A1" \ar@{-}@/^{4ex}/|{\scriptscriptstyle M}   "A5"
\POS"A3" \ar@{-}@/_{4ex}/|{\scriptscriptstyle N\ }   "A7"
\POS"A5" \ar@{-}|{\scriptscriptstyle L_1}   "A7"
\POS"A1" \ar@{-}@/^{10ex}/|{\scriptscriptstyle \,L_2\,}   "A3"
\POS"A1" \ar@{-}|{\scriptscriptstyle L'_1}   "A7"
\POS"A3" \ar@{-}|{\scriptscriptstyle L'_2}   "A5"
\endxy
&\qquad\qquad&
\xy/r4pc/:{\xypolygon8"A"{~<{}~>{-}{\scriptstyle\bullet}}}
*+{{\scriptstyle \bullet}},
\POS"A6" \ar@{-}@/^{14ex}/|(0.9){\scriptscriptstyle \ M}   "A7"
\POS"A6" \ar@{-}@/_{4ex}/|(0.7){\scriptscriptstyle N}   "A2"
\POS"A2" \ar@{-}@/^{4ex}/|{\scriptscriptstyle L_1}   "A7"
\POS"A6" \ar@{-}@/^{0.5ex}/|(0.7){\scriptscriptstyle L_2\,}   "A0"
\POS"A6" \ar|-(0.75){\SelectTips{cm}{}\object@{+}}@{-}@/_{1ex}/|(0.5){\scriptscriptstyle \  L_3}   "A0"
\POS"A2" \ar@{-}@/_{5ex}/|{\scriptscriptstyle \,L'_1}   "A6"
\POS"A6" \ar@{-}_{\scriptscriptstyle \,L'_2=0}   "A7"
\endxy
\end{array} 
\] 
\caption{Two examples of exchange pairs}\label{fig exchange}
\end{figure}   
\end{subsection}

\begin{subsection}{Cluster-tilted algebra}
Let $T$ be any triangulation of the punctured polygon. Then the
endomorphism algebra $\End_{\mathcal{C}}(T)^{op}$ is called {\em
  cluster-tilted algebra}. 
By a result of \cite{BMR1}, the functor $\Hom_{\mathcal{C}}(\tau^{-1}\,T,-)$
induces an equivalence of 
categories $\varphi_T:\mathcal{C}/T \to \textup{mod}\,
\End_{\mathcal{C}}(T)^{op} $.
Labeling the edges in $T$ by $T_1,T_2,\ldots,T_n$, we get the dimension
vector $\underline{\dim}\,\varphi_T(M_{a,b}^\ze)$ of a module
$\varphi_T(M_{a,b}^\ze)$ by 
\begin{equation} \label{eins}
  (\underline{\dim}\,\varphi_T(M_{a,b}^\ze))_i =
\dim\Hom_{\mathcal{C}}(\tau^{-1}\,T_i, M_{a,b}^\ze) =
\dim\Ext^1_{\mathcal{C}}(M_{a,b}^\ze,T_i) = e(M_ {a,b}^\ze,T_i),
\end{equation}
where the last equation holds by Theorem \ref{thm 2}.

We illustrate this in an example.
Let $n=4$ and $T$ be the triangulation 
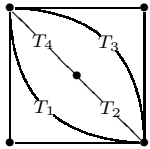
\[ \xy/r3pc/:{\xypolygon4"A"{~<{}~>{-}{\scriptstyle\bullet}}}
*+{{\scriptstyle \bullet}},
\POS"A2" \ar@{-}@/_{4ex}/|{T_1}   "A4"
\POS"A2" \ar@{-}@/^{4ex}/|{T_3}   "A4"
\POS"A2" \ar@{-}|{T_4}   "A0"
\POS"A4" \ar@{-}|{T_2}   "A0"
\endxy\] 
Then the cluster-tilted algebra $\End_\mathcal{C}(T)^{op}$ is the
quotient of  the path algebra 
of the quiver 
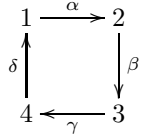
\[\xymatrix{1\ar[r]^\za&2\ar[d]^\zb\\
4\ar[u]^\zd&3\ar[l]^\zg
}\]
by the ideal generated by the paths $\za\,\zb\,\zg$, $\zb\,\zg\,\zd$,
$\zg\,\zd\,\za$ and $\zd\,\za\,\zb$.
The Auslander-Reiten quiver of the category $\mathcal{C}$ is
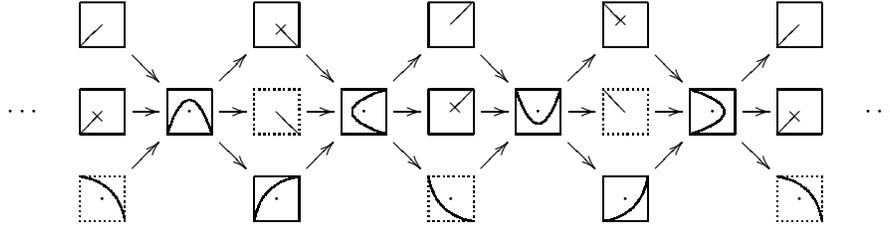
\[
\xymatrix@R=10pt@C=10pt{ 
&\xy/r1pc/:{\xypolygon4"A"{~<{}~>{-}{}}}
\POS"A3" \ar@{-}   "A0"
\endxy
\ar[rd]&&
\xy/r1pc/:{\xypolygon4"A"{~<{}~>{-}{}}}
\POS"A4" \ar|-(0.8){\SelectTips{cm}{}\object@{+}}@{-}    "A0"
\endxy
\ar[rd]&&
\xy/r1pc/:{\xypolygon4"A"{~<{}~>{-}{}}}
\POS"A1" \ar@{-}   "A0"
\endxy
\ar[rd]&&
\xy/r1pc/:{\xypolygon4"A"{~<{}~>{-}{}}}
\POS"A2" \ar|-(0.8){\SelectTips{cm}{}\object@{+}}@{-}   "A0"
\endxy
\ar[rd]
&&\xy/r1pc/:{\xypolygon4"A"{~<{}~>{-}{}}}
\POS"A3" \ar@{-}   "A0"
\endxy
\\
\cdots&
\xy/r1pc/:{\xypolygon4"A"{~<{}~>{-}{}}}
\POS"A3" \ar|-(0.8){\SelectTips{cm}{}\object@{+}}@{-}   "A0"
\endxy
\ar[r]
&\xy/r1pc/:{\xypolygon4"A"{~<{}~>{-}{}}}
*+{\cdot},
\POS"A3" \ar@{-}@/^{3ex}/   "A4"
\endxy
\ar[rd]\ar[ru]\ar[r]&
\xy/r1pc/:{\xypolygon4"A"{~<{}~>{.}{}}}
\POS"A4"  \ar@{-}    "A0"
\endxy
\ar[r]
&\xy/r1pc/:{\xypolygon4"A"{~<{}~>{-}{}}}
*+{\cdot},
\POS"A4" \ar@{-}@/^{3ex}/   "A1"
\endxy
\ar[rd]\ar[ru]\ar[r]&
\xy/r1pc/:{\xypolygon4"A"{~<{}~>{-}{}}}
\POS"A1" \ar|-(0.8){\SelectTips{cm}{}\object@{+}}@{-}   "A0"
\endxy
\ar[r]
&\xy/r1pc/:{\xypolygon4"A"{~<{}~>{-}{}}}
*+{\cdot},
\POS"A1" \ar@{-}@/^{3ex}/   "A2"
\endxy
\ar[rd]\ar[ru]\ar[r]&
\xy/r1pc/:{\xypolygon4"A"{~<{}~>{.}{}}}
\POS"A2"  \ar@{-}    "A0"
\endxy
\ar[r]
&\xy/r1pc/:{\xypolygon4"A"{~<{}~>{-}{}}}
*+{\cdot},
\POS"A2" \ar@{-}@/^{3ex}/   "A3"
\endxy
\ar[rd]\ar[ru]\ar[r]&
\xy/r1pc/:{\xypolygon4"A"{~<{}~>{-}{}}}
\POS"A3" \ar|-(0.8){\SelectTips{cm}{}\object@{+}}@{-}   "A0"
\endxy
&\cdots
\\
&\xy/r1pc/:{\xypolygon4"A"{~<{}~>{.}{}}}
*+{\cdot},
\POS"A2" \ar@{-}@/^{1ex}/   "A4"
\endxy
\ar[ru]&&
\xy/r1pc/:{\xypolygon4"A"{~<{}~>{-}{}}}
*+{\cdot},
\POS"A3" \ar@{-}@/^{1ex}/   "A1"
\endxy
\ar[ru]&&
\xy/r1pc/:{\xypolygon4"A"{~<{}~>{.}{}}}
*+{\cdot},
\POS"A4" \ar@{-}@/^{1ex}/   "A2"
\endxy
\ar[ru]&&
\xy/r1pc/:{\xypolygon4"A"{~<{}~>{-}{}}}
*+{\cdot},
\POS"A1" \ar@{-}@/^{1ex}/   "A3"
\endxy
\ar[ru]&&
\xy/r1pc/:{\xypolygon4"A"{~<{}~>{.}{}}}
*+{\cdot},
\POS"A2" \ar@{-}@/^{1ex}/   "A4"
\endxy
}
\]

For the four tagged edges
of the triangulation $T$, we have drawn the borders of the
punctured polygons in the above picture
as  dotted lines. Deleting these positions and using equation
(\ref{eins}), we obtain the Auslander-Reiten 
quiver of $\End_{\mathcal{C}}(T)^{op}$

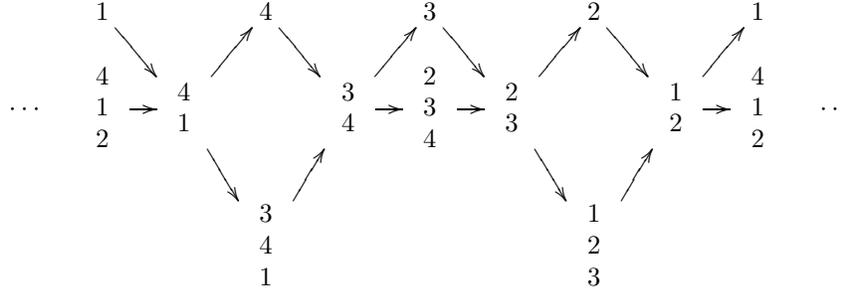
\[
\xymatrix@R=10pt@C=10pt{ 
&
1\ar[rd]&&
4\ar[rd]&&
3\ar[rd]&&
2\ar[rd]&&
1
\\
\cdots&
{\begin{array}{c} 
4\\1\\2
\end{array}  }
\ar[r]&
{\begin{array}{c} 4\\1
\end{array}  }
\ar[ru]\ar[rd]&&
{\begin{array}{c} 3\\4
\end{array}  }
\ar[ru]\ar[r]&
{\begin{array}{c} 
2\\3\\4
\end{array}  }
\ar[r]&
{\begin{array}{c} 2\\3
\end{array}  }
\ar[ru]\ar[rd]&&
{\begin{array}{c} 1\\2
\end{array}  }
\ar[ru]\ar[r]&
{\begin{array}{c} 
4\\1\\2
\end{array}  }
&\cdots
\\
&&&
{\begin{array}{c} 
3\\4\\1
\end{array}  }
\ar[ru]&&
&&
{\begin{array}{c} 
1\\2\\3
\end{array}  }
\ar[ru]&&
}
\]
where modules are represented by their Loewy series.

\end{subsection} 
\end{section} 
{} 

\bigskip 
\bigskip 

\noindent 
Department of  Mathematics and    Statistics\\
University of   Massachusetts at Amherst\\
Amherst, MA 01003--9305, USA\\ 
E-mail address: {\tt schiffler@math.umass.edu}

\end{document}